\documentclass[11pt]{article}
\usepackage{bm,amssymb,amsmath}

\newtheorem{definition}{Definition}[section]
\newtheorem{theorem}{Theorem}[section]
\newtheorem{proposition}[theorem]{Proposition}
\newtheorem{lemma}[theorem]{Lemma}
\newtheorem{corollary}[theorem]{Corollary}
\newtheorem{example}{Example}[section]

\newcommand{\point}[1]{\bm{#1}}
\newcommand{\p}[1]{\point{#1}}
\newcommand{\field}[1]{\mathcal{#1}}
\newcommand{\pset}[1]{\mathbb{#1}}
\newcommand{\ideal}[1]{\mathcal{#1}}
\newcommand{\bases}[1]{\langle #1 \rangle}
\newcommand{\fk}{\field{K}}
\newcommand{\res}{{\rm res}}
\newcommand{\lpp}{{\rm lm}}
\newcommand{\xx}{x_1,\ldots,x_n}
\newcommand{\nf}{{\rm nform}}
\DeclareMathOperator{\cls}{cls}
\DeclareMathOperator{\coef}{coef}
\DeclareMathOperator{\ini}{ini}
\DeclareMathOperator{\lc}{lc}
\DeclareMathOperator{\lt}{lt}

\DeclareMathOperator{\lv}{lv}
\DeclareMathOperator{\pquo}{pquo}
\DeclareMathOperator{\prem}{prem}
\DeclareMathOperator{\sat}{sat}
\DeclareMathOperator{\zero}{Zero}
\DeclareMathOperator{\nform}{nform}

\begin{document}

\title{On the Connection Between Ritt Characteristic Sets and Buchberger-Gr\"obner Bases}

\author{Dongming Wang \\[5pt] \normalsize
School of Mathematics and Systems Science, \\ \normalsize
Beihang University, Beijing 100191, China \\[3pt] \normalsize
Centre National de la Recherche Scientifique, \\ \normalsize
75794 Paris cedex 16, France
}
\date{}

\maketitle

\begin{abstract}
For any polynomial ideal $\ideal{I}$,
let the minimal triangular set contained in the reduced Buchberger-Gr\"obner basis of $\ideal{I}$ with respect to the purely lexicographical term order be called the W-characteristic set of $\ideal{I}$. In this paper, we establish a strong connection between Ritt's characteristic sets and Buchberger's Gr\"obner bases of polynomial ideals by showing that the W-characteristic set $\pset{C}$ of $\ideal{I}$ is a Ritt characteristic set of $\ideal{I}$ whenever $\pset{C}$ is an ascending set, and a Ritt characteristic set of $\ideal{I}$ can always be computed from $\pset{C}$ with simple pseudo-division when $\pset{C}$ is regular. We also prove that under certain variable ordering, either the W-characteristic set of $\ideal{I}$ is normal, or irregularity occurs for the $j$th, but not the $(j+1)$th, elimination ideal of $\ideal{I}$ for some $j$. In the latter case, we provide explicit pseudo-divisibility relations, which lead to nontrivial factorizations of certain polynomials in the Buchberger-Gr\"obner basis and thus reveal the structure of such polynomials. The pseudo-divisibility relations may be used to devise an algorithm to
decompose arbitrary polynomial sets into normal triangular sets based on Buchberger-Gr\"obner bases computation.

\end{abstract}

\section{Introduction}\label{sec:intro}

In his classical works \cite{ri32,ri50}, Ritt introduced the concept of characteristic sets for polynomial and differential polynomial ideals. This concept plays a central role in Ritt's theory of differential algebra. Ritt did not provide any effective way for the construction of
characteristic sets even for non-prime polynomial ideals. Wu \cite{wu94} developed Ritt's theory and method by working out efficient algorithms for the computation of characteristic sets of (differential) polynomial sets instead of (differential) polynomial ideals. On the other hand, in his Ph.D.\ thesis \cite{bu65} Buchberger introduced the concept of Gr\"obner bases for polynomial ideals and proposed an algorithm for effective computation of Gr\"obner bases. Over the past three decades the theories and methods of characteristic sets and Gr\"obner bases have been studied extensively and independently by many researchers. They have become fundamental tools for computational commutative algebra and algebraic geometry (see, e.g., \cite{bw93,clo96,ei95,ms04}).

Characteristic sets and Gr\"obner bases are rather different in terms of structure and properties. The algorithms for their computation are distinguished by their elimination strategies (elimination of variables vs.\ elimination of terms) and reduction steps (pseudo-division vs.\ Buchberger reduction). It was not known what inherent connections may exist between the characteristic sets and the Gr\"obner bases of arbitrary polynomial ideals. For any polynomial ideal $\ideal{I}$, we call the minimal triangular set contained in the reduced Buchberger-Gr\"obner basis of $\ideal{I}$ with respect to the purely lexicographical term order the \emph{W-characteristic set} of $\ideal{I}$. This paper establishes, for the first time, a strong connection between Ritt's characteristic sets and Buchberger's Gr\"obner bases of polynomial ideals by showing that the W-characteristic set $\pset{C}$ of $\ideal{I}$ is a Ritt characteristic set of $\ideal{I}$ whenever $\pset{C}$ is an ascending set, and a Ritt characteristic set of $\ideal{I}$ can always be computed from $\pset{C}$ with simple pseudo-division when $\pset{C}$ is regular. We also prove that under certain variable ordering, either the W-characteristic set of $\ideal{I}$ is normal, or irregularity occurs for the $j$th, but not the $(j+1)$th, elimination ideal of $\ideal{I}$ for some $j$. In the latter case, we provide explicit pseudo-divisibility relations, which lead to nontrivial factorizations of certain polynomials in the Buchberger-Gr\"obner basis and thus reveal the structure of such polynomials. It is not clear if the uniquely defined W-characteristic set of $\ideal{I}$ contains all what is needed for the construction of a Ritt characteristic set of $\ideal{I}$,
but it does appear to contain sufficient characteristic information about $\ideal{I}$, in the abnormal case. The pseudo-divisibility relations may be used to devise an algorithm to
decompose arbitrary polynomial sets into normal triangular sets based on Buchberger-Gr\"obner bases computation.

The strong connection between the Ritt characteristic set and the Buchberger-Gr\"obner basis of a polynomial ideal shown in this paper is simple, yet deep and surprising. It creates a route for the study of characteristic sets of polynomial ideals using the theory and method of Buchberger-Gr\"obner bases, and vice versa. Our work also illustrates how remarkable results in polynomial ideal theory can be obtained from the interplay of the two conceptually and operationally different methods. Further investigations on the W-characteristic sets of polynomial ideals are likely to
help deepen our understanding of the structural properties of both characteristic sets and Buchberger-Gr\"obner bases.

\section{Preliminaries}\label{sec:pre}

Let $\fk$ be a field and $x_1, \ldots, x_n$ be $n$ variables
with a fixed order $x_1 <_{{\rm plex}} \cdots <_{{\rm plex}} x_n$. Denote by $\fk[x_1, \ldots, x_n]$
the ring of polynomials in $x_1, \ldots, x_n$ with coefficients in $\fk$.

Let $F\in\fk[x_1, \ldots, x_n]$ be any nonzero polynomial. Denote by
$\deg(F, x_k)$ the \emph{degree} of $F$ in $x_k$ and by $\lc(F, x_k)$
the \emph{leading coefficient} of $F$ with respect to $x_k$.
Let $m=\deg(F,x_k)$ and $G$ be any other polynomial
of degree $l$ in $x_k$. Pseudo-dividing $G$ by $F$,
considered as polynomials in $x_k$, one can obtain
two polynomials $Q$ and $R$ in $\fk[x_1, \ldots, x_n]$ such that
\begin{equation}\label{rmdr}
I^qG=Q F+R,
\end{equation}
where
\[\begin{array}{l}\smallskip
I=\lc(F,x_k), ~~q=\max(l-m+1,0), \\
\deg(R,x_k)<m, ~~\deg(Q,x_k)=\max(l-m, -1).
\end{array}\]
In case $m=0$, $R=0$ and $Q=G^lF$. The uniquely determined polynomials $Q$ and $R$ are called the
\emph{pseudo-quotient} and the \emph{pseudo-remainder} of $G$ with respect to $F$ in $x_k$,
denoted by $\pquo(G, F, x_k)$ and $\prem(G, F, x_k)$, respectively.

The polynomial $G$ is said to be \emph{R-reduced} with respect to $F$ in $x_k$ if $l<m$.
When $G$ is R-reduced with respect to $F$ in $x_k$, $R=G$ and $Q=0$.

Let $P\in \fk[x_1, \ldots, x_n]\setminus\fk$ be any nonconstant polynomial. The biggest index $p$ such that
$x_p$ actually occurs in $P$ is called the \emph{class} of $P$,
denoted by $\cls(P)$, and the variable $x_p$ is called the \emph{leading
variable} of $P$, denoted by $\lv(P)$. The class of any constant polynomial in $\fk$
is defined to be $0$. Let $\cls(P) = p >0$; then $P$ can be
written as $P = I x_p^d + H$ with $I \in \fk[x_1,\ldots,x_{p-1}]$
and $\deg(H, x_k)< d=\deg(P, x_p)$. The leading coefficient $I$ of $P$ with respect to $x_p$ is called the \emph{initial} of $P$, denoted by $\ini(P)$.

When a polynomial $G$ is R-reduced with respect to $P$ in $x_p=\lv(P)$, we simply say that $G$ is R-reduced with respect to $P$ (without mentioning $x_p$).

\subsection{Triangular sets and characteristic  sets}

 To represent an ordered set, we enclose its elements using a pair of square brackets instead of braces.

\begin{definition}\em
A finite nonempty ordered set
$[T_1, \ldots, T_r]$
 of nonconstant polynomials in $\fk[x_1, \ldots, x_n]$ is called a \emph{triangular set} if
$0<\cls(T_1)<\cdots<\cls(T_r)$.

A finite nonempty ordered set
$\pset{A}=[A_1, \ldots, A_r]$
 of nonzero polynomials in $\fk[x_1, \ldots, x_n]$ is called an \emph{ascending set} if either $r=1$ and $A_1\in\fk$, or $\pset{A}$ is a triangular set and $A_j$ is R-reduced with respect to $A_i$ for every pair $i<j$ and $j=2,\ldots,r$.
\end{definition}

Let $\pset{T}=[T_1, \ldots, T_r]$ be any triangular set
and $P$ be any polynomial in $\fk[x_1, \ldots, x_n]$. $P$ is said to be \emph{R-reduced}
 with respect to $\pset{T}$ if $P$ is R-reduced with respect to every $T_i\in\pset{T}$, i.e.,
$\deg(P, \lv(T_i))<\deg(T_i, \lv(T_i))$ for all $1\leq i\leq r$.

The polynomial
\[R=\prem(\cdots\prem(P,T_r,\lv(T_r)), \ldots, T_1, \lv(T_1)),\]
denoted simply by $\prem(P, \pset{T})$, is called
the \emph{pseudo-remainder} of $P$ with respect to $\pset{T}$.
From (\ref{rmdr}), one can easily derive the following
\emph{pseudo-remainder formula}
\begin{equation}\label{rmdr:TS}
I_1^{q_1}\cdots I_r^{q_r} P=Q_1T_1+\cdots + Q_r T_r+R,
\end{equation}
where each $q_i$ is a nonnegative integer and
\[I_i=\ini(T_i), ~~Q_i\in\fk[x_1, \ldots, x_n], \quad i=1, \ldots, r.\]

\begin{definition}\em
For any two nonzero polynomials $F$ and $G$ in $\fk[x_1, \ldots, x_n]$, $F$ is said
to have \emph{lower rank} than $G$, denoted as
$F\prec G$ or $G\succ F$,
if either $\cls(F)<\cls(G)$, or $\cls(F)=\cls(G)>0$ and $\deg(F,\lv(F))<\deg(G,\lv(G))$.
In this case, $G$ is said to have \emph{higher rank} than $F$.

If neither $F\prec G$ nor $G\prec F$, then $F$ and $G$ are said to have
the \emph{same rank}, denoted as $F\sim G$.
\end{definition}

\begin{definition}\em
For any two ascending sets
\[\pset{A}=[A_1, \ldots, A_r], ~~\pset{A}'=[A_1', \ldots, A_{r'}'],\]
$\pset{A}$ is said to have \emph{higher rank} than
$\pset{A}'$, denoted as $\pset{A}\succ\pset{A}'$ or $\pset{A}'\prec\pset{A}$,
if one of the following conditions holds:
\begin{enumerate}
\item[(a)]  there exists an integer $j\leq \min(r, r')$ such that
\[A_1\sim A_1', \ldots, A_{j-1}\sim A_{j-1}',~~\mbox{while}~~A_j\succ A_j';\]
\item[(b)]  $r'>r$ and
$A_1\sim A_1', \ldots, A_r\sim A_r'$.
\end{enumerate}
In this case, $\pset{A}'$ is said to have \emph{lower rank} than $\pset{A}$.
If neither $\pset{A}\prec\pset{A}'$ nor $\pset{A}'\prec\pset{A}$, then
$\pset{A}$ and $\pset{A}'$ are said to have the \emph{same rank}, denoted as
$\pset{A}\sim \pset{A}'$. In this case,
\[r=r',~~\mbox{and}~~A_1\sim A_1', \ldots, A_r\sim A_r'.\]
\end{definition}

For any polynomial set $\pset{P}=\{P_1,\ldots,P_s\}$ in $\fk[x_1,\ldots,x_n]$, denote by
$\langle \pset{P}\rangle$ or $\langle P_1,\ldots,P_s\rangle$ the ideal generated by the
polynomials $P_1,\ldots,P_s$ in $\fk[x_1,\ldots,x_n]$ and by $\zero(\pset{P})$ the set
of all common zeros (in some extension field of $\fk$) of $P_1,\ldots,P_s$.

\begin{definition}\em
Let $\pset{P}$ be any finite nonempty polynomial set in $\fk[x_1,\ldots,x_n]$.
With respect to the rank $\succ$, any minimal ascending set contained in $\langle \pset{P}\rangle$ is called
a \emph{Ritt characteristic set} of $\langle \pset{P}\rangle$.
\end{definition}

\begin{lemma}\label{rittlem}
For any finite nonempty polynomial set $\pset{P}$ in $\fk[x_1,\ldots,x_n]$,
an ascending set $\pset{A}$ contained in $\langle \pset{P}\rangle$
is a \emph{Ritt characteristic set} of $\langle \pset{P}\rangle$ if and only
if $\prem(P, \pset{A})\equiv 0$ for all
$P\in\langle \pset{P}\rangle$.
\end{lemma}
\emph{Proof.} See \cite[Theorem 5.3.3]{mi93}. $\square$

\begin{definition}\em\label{DEF:sat}
For any triangular set $\pset{T}=[T_1, \ldots, T_r]$ in $\fk[x_1, \ldots, x_n]$ with $I_i=\ini(T_i)$ for $1\leq i\leq r$, the \emph{saturated ideal}
of $\pset{T}$ is defined to be
\[\sat(\pset{T})=\{F\in\fk[x_1, \ldots, x_n]\,|
~(I_1\cdots I_r)^qF\in\bases{\pset{T}}~\mbox{for~some~integer}~q>0\}.\]
\end{definition}

\subsection{Regular sets and normal triangular sets}

For two polynomials $F$ and $G$ in $\fk[x_1, \ldots, x_n]$, denote by
$\res(F,G,x_k)$ the \emph{resultant} of $F$ and $G$ with respect to $x_k$. For any polynomial $F$
and triangular set $\pset{T}=[T_1, \ldots, T_r]$, define
$\res(F, \pset{T})=\res(\cdots\res(F,T_r,\lv(T_r)),\ldots,T_1,\lv(T_1))$. It is easy to see that there exist polynomials $A, B_1,\ldots,B_r\in\fk[x_1, \ldots, x_n]$ such that
\begin{equation}\label{resform}
\res(F, \pset{T})=AF+B_1T_1+\cdots +B_rT_r
\end{equation}
(cf.\ \cite[Lemmas 7.2.1 and 7.2.2]{mi93} and \cite[Lemma 4.3.2]{wa01}), so $\res(F, \pset{T})\in\langle F, T_1, \ldots, T_r\rangle$.

\begin{definition}\label{regdef}\em
  Let $\pset{T}=[T_1, \ldots, T_r]$ be any triangular set in $\fk[x_1, \ldots, x_n]$.
  $\pset{T}$ is said to be \emph{regular} or called a \emph{regular set} if $\res(\ini(T_j), [T_1,\ldots,T_{j-1}])\not\equiv 0$ for all $j=2, \ldots, r$.
  $\pset{T}$ is said to be \emph{normal} if $\deg(\ini(T_j), \lv(T_i)) = 0$
  for all $i<j$ and $j=2, \ldots, r$.
\end{definition}

Regular sets, also known as \emph{regular chains}, have other equivalent definitions (see, e.g.,
\cite[p.\,114]{wa01}). Here they are defined by means of resultants for the convenience of proof in Section~3. A triangular set $\pset{T}$ as in Definition~\ref{regdef} is regular
if and only if the image of $\ini(T_i)$ is neither zero nor a zero divisor in the quotient ring
$\fk[x_1, \ldots, x_n]/\sat([T_1,\ldots,T_{i-1}])$ for all $i=2, \ldots, r$.

\begin{proposition}\label{regpro}
Every normal triangular set in $\fk[x_1, \ldots, x_n]$ is regular.
\end{proposition}
\emph{Proof.}
It is obvious.
$\square$

\begin{proposition}\label{satpro}
  Let $\pset{T}$ be any triangular set in $\fk[x_1, \ldots, x_n]$.
  Then $\pset{T}$ is regular if and only if $\prem(P, \pset{T})\equiv 0$ for all $P\in\sat(\pset{T})$.
\end{proposition}
\emph{Proof.}
See \cite[Theorem~6.2.4]{wa01} and references therein.
$\square$

\begin{lemma}\label{satlem}
  Let $\pset{T}=[T_1, \ldots, T_r]$ with $p_r=\cls(T_r)<n$ be any regular set and $P=P_dx_m^d+\cdots+P_1x_m+P_0$ with $m>p_r$ and $d=\deg(P,x_m)>0$ be any polynomial in $\fk[x_1, \ldots, x_n]$.
  Then $\prem(P, \pset{T})\equiv 0$ if and only if $\prem(P_j, \pset{T})\equiv 0$ for all $0\leq j\leq d$.
\end{lemma}
\emph{Proof.} Let $I_i=\ini(T_i)$ for $1\leq i\leq r$.
Suppose that $\prem(P, \pset{T})\equiv 0$. Then according to \eqref{rmdr:TS} there exist a power product $K$ of $I_1, \ldots, I_r$ (i.e., $K=I_1^{q_1}\cdots I_r^{q_r}$ for some nonnegative integers $q_1, \ldots, q_r$) and polynomials $Q_1,\ldots,Q_r$ in $\fk[x_1, \ldots, x_n]$ such that $KP=
Q_1T_1+\cdots + Q_rT_r$. It follows that \[KP_j=\coef(Q_1T_1+\cdots + Q_rT_r, x_m^j)=\coef(Q_1,x_m^j)\,T_1+\cdots + \coef(Q_r,x_m^j)\,T_r\in\langle\pset{T}\rangle\]
for $0\leq j\leq d$, where $\coef(F, x_m^j)$ denotes the coefficient of $F$ in $x_m^j$. Therefore, $P_j\in\sat(\pset{T})$ by definition. As $\pset{T}$ is regular, we have $\prem(P_j, \pset{T})\equiv 0$ according to Proposition~\ref{satpro}.

To show the other direction, suppose that $\prem(P_j, \pset{T})\equiv 0$ for all $0\leq j\leq d$. Then $P_j\in\sat(\pset{T})$ for all $j$. This implies that $P\in\sat(\pset{T})$. Since $\pset{T}$ is regular, we have
$\prem(P, \pset{T})\equiv 0$, again by Proposition~\ref{satpro}. The proof is complete.
$\square$

\begin{corollary}\label{satcor}
  Let $\pset{T}$ be any regular set and $P$ and $F$ be any two polynomials in $\fk[x_1, \ldots, x_n]$.
   If $\prem(P, \pset{T})=\prem(F, \pset{T})\equiv 0$, then $\prem(P+F, \pset{T})\equiv 0$.
\end{corollary}
\emph{Proof.} Let $y$ be a variable not occurring in $P, F$, and $\pset{T}$. By Lemma~\ref{satlem},
$\prem(P, \pset{T})=\prem(F, \pset{T})\equiv 0$ implies that $\prem(Py+F, \pset{T})\equiv 0$. The corollary is proved by taking $y=1$. $\square$

\medskip
Lemma~\ref{satlem} and Corollary~\ref{satcor} do not hold when $\pset{T}$ is not regular. This can be seen from  the simple example, where $\pset{T}=[x_1^2, x_1x_2]$, $P=x_3-x_2^2$, and $F=x_2^2$; it is easy to see that $\prem(P, \pset{T})=\prem(F, \pset{T})\equiv 0$, but $\prem(\lc(P, x_3), \pset{T})\not\equiv 0$ and $\prem(P+F, \pset{T})\not\equiv 0$.

\subsection{Buchberger-Gr\"obner bases}

Two distinct monomials $x_1^{i_1}\cdots x_n^{i_n}$ and $x_1^{j_1}\cdots x_n^{j_n}$ in $x_1,\ldots,x_n$ are ordered as
\[x_1^{i_1}\cdots x_n^{i_n} <_{{\rm plex}} x_1^{j_1}\cdots x_n^{j_n} \quad\mbox{or}\quad
x_1^{j_1}\cdots x_n^{j_n} >_{{\rm plex}}   x_1^{i_1}\cdots x_n^{i_n},\]
if there exists an integer $k$ ($1\leq k\leq n$) such that
\[i_n=j_n, \ldots , i_{k+1}=j_{k+1} \quad\mbox{while}\quad i_k<j_k.\]
Under $<_{{\rm plex}}$, all the monomials in $x_1,\ldots,x_n$ are ordered,
and so are the terms of any nonzero polynomial in $\fk[x_1,\ldots,x_n]$.
We call $<_{{\rm plex}}$ the \emph{purely lexicographical order} (plex)
of monomials or terms.

Any nonzero polynomial $P$ in $\fk[x_1,\ldots,x_n]$ can be written
in the form
\[P=\sum_{l=1}^t a_lx_1^{i_{l1}}\cdots x_n^{i_{ln}}
\] with
\[\begin{array}{l}\smallskip
a_1\neq 0, \ldots , a_t\neq 0, ~~a_i\in\fk, \\
x_1^{i_{11}}\cdots x_n^{i_{1n}}>_{{\rm plex}}  \cdots>_{{\rm plex}}   x_1^{i_{t1}}\cdots x_n^{i_{tn}}.
\end{array}\]
We call $x_1^{i_{11}}\cdots x_n^{i_{1n}}$ the \emph{leading monomial},
$a_1 x_1^{i_{11}}\cdots x_n^{i_{1n}}$ the \emph{leading term}, $a_1$ the \emph{leading
coefficient} of $P$, and $a_l$ the \emph{coefficient} of $P$ in $x_1^{i_{l1}}\cdots x_n^{i_{ln}}$, denoted by $\lpp(P)$,
$\lt(P)$, $\lc(P)$, and $\coef(P, x_1^{i_{l1}}\cdots x_n^{i_{ln}})$, respectively.
If $Q$ is another nonzero polynomial in $\fk[x_1,\ldots,x_n]$, we order $P$ and $Q$ as
$P<_{{\rm plex}} Q$ or $Q>_{{\rm plex}} P$ if $\lpp(P)<_{{\rm plex}} \lpp(Q)$.

\begin{definition}\em
Let $\pset{P}$ be any finite nonempty polynomial set and $G$ be any polynomial in $\fk[\xx]$.
$G$ is said to be \emph{B-reducible} with respect to $\pset{P}$ if there exist a polynomial
$P\in\pset{P}$ and a monomial $\lambda$ such that
$\coef(G, \lambda\,\lpp(P))\neq 0$.
If no such $P$ and $\lambda$ exist, $G$ is said to be \emph{B-reduced}
or in \emph{normal form} with respect to $\pset{P}$.
\end{definition}

If $G$ is B-reducible with respect to $\pset{P}$, then one can find a polynomial $P\in\pset{P}$
with the monomial $\lambda\,\lpp(P)$ maximal (with respect to the term order $<_{{\rm plex}}$)
such that
\[G=b\, \lambda \, P+H,\]
where
\[b=\frac{\coef(G, \lambda\,\lpp(P))}{\lc(P)}.\]

If $H$ is B-reducible with respect to $\pset{P}$, then one can
reduce $H$ to another polynomial in the same way by choosing $P, b$,
and $\lambda$. Such a
process will terminate after a finite number of reduction
steps. The finally obtained polynomial $N$ will be B-reduced with respect to $\pset{P}$.
In this case, one gets a formula of the form
\[G=Q_1P_1+\cdots+Q_sP_s+N,\]
in which $P_j\in\pset{P}$, $Q_j, N\in\fk[\xx]$ and $N$ is B-reduced with respect to $\pset{P}$.
The polynomial $N$ is called the \emph{normal form}
of $G$ with respect to $\pset{P}$ and denoted by $\nf(G,\pset{P})$.

\begin{definition}\em
Let $\pset{P}$ be an arbitrary finite and nonempty set of nonzero polynomials
in $\fk[x_1,\ldots,x_n]$. A polynomial set $\pset{G}$ in $\fk[x_1,\ldots,x_n]$ is called
the \emph{reduced Buchberger-Gr\"obner basis} of $\langle\pset{P}\rangle$ or $\pset{P}$ with respect to the plex term order determined by $x_1<_{{\rm plex}}\cdots <_{{\rm plex}} x_n$, if
\begin{enumerate}
\item[{\rm (a)}] for all $P\in\fk[x_1,\ldots,x_n]$, $P\in\langle\pset{P}\rangle$ if and only if $\nf(P, \pset{G})=0$;
\item[{\rm (b)}] every polynomial $G\in\pset{G}$ is monic and B-reduced with respect to
$\pset{G}\setminus \{G\}$;
\item[{\rm (c)}] $\langle\pset{G}\rangle=\langle\pset{P}\rangle$.
\end{enumerate}
\end{definition}

The reduced Buchberger-Gr\"obner basis of $\langle\pset{P}\rangle$ is unique and can be computed from $\pset{P}$ by using Buchberger's algorithm \cite{bu85}. What is called Buchberger-Gr\"obner basis here was named Gr\"obner basis by Buchberger after his Ph.D.\ advisor Wolfgang Gr\"obner. The author feels that the basis should be named more appropriately also after Bruno Buchberger for his outstanding contributions to the development of the theory and method of Gr\"obner bases.

\section{Main Results}\label{sec:main}

For any (finite or infinite) polynomial set $\pset{F}\subset\fk[x_1,\ldots,x_n]$ and
$0\leq j\leq k\leq n$, let $\pset{F}^{\langle j,\ldots,k\rangle}$ stand for
$(\pset{F}\cap\fk[x_1,\ldots,x_k])\setminus (\pset{F}\cap\fk[x_1,\ldots,x_{j-1}])$.
When $j=k$, $\pset{F}^{\langle j,\ldots,k\rangle}$ is written as $\pset{F}^{\langle k\rangle}$.

\begin{definition}\label{wcharsetdef}\em
Let $\pset{P}$ be an arbitrary finite and nonempty set of nonzero polynomials
in $\fk[x_1,\ldots,x_n]$ and $\pset{G}$ be the reduced Buchberger-Gr\"obner basis of $\pset{P}$ with respect to the plex term order determined by $x_1<_{{\rm plex}}\cdots <_{{\rm plex}} x_n$.
The set
\[\bigcup_{i=0}^n\left\{G\,\big|\big.~G\in\pset{G}^{\langle i\rangle};\, G'>_{{\rm plex}} G, ~\mbox{for all}~ G'\in\pset{G}^{\langle i\rangle}\setminus \{G\}\right\}\]
of polynomials, ordered by $<_{{\rm plex}}$, is called the \emph{W-characteristic set} of $\langle\pset{P}\rangle$.
\end{definition}

The above-defined W-characteristic set $\pset{C}$ is obviously a triangular set, but it is not necessarily an ascending set. $\pset{C}$ is minimal in the sense that (i) each element $C$ of $\pset{C}$ has the lowest plex order among all those polynomials in the Buchberger-Gr\"obner basis $\pset{G}$ which have the same leading variable as $C$ and (ii) the number of elements in $\pset{C}$ is the maximum possible.
Owing to the uniqueness of the reduced Buchberger-Gr\"obner basis, the W-characteristic set of a polynomial ideal is uniquely defined.
We will see that W-characteristic sets can be effectively used to bridge Ritt characteristic sets and Buchberger-Gr\"obner bases.

\subsection{Construction of Ritt characteristic sets}

In what follows, let $\pset{P}$ be an arbitrary finite and nonempty set of nonzero polynomials
in $\fk[x_1,\ldots,x_n]$, let $\pset{G}$ be the reduced Buchberger-Gr\"obner basis of $\pset{P}$ with respect to the plex term order determined by $x_1<_{{\rm plex}}\cdots <_{{\rm plex}} x_n$, and assume that $\pset{G}\neq [1]$, so $\pset{C}\neq [1]$, whenever needed. The following proposition shows that the W-characteristic set of $\langle\pset{P}\rangle$ possesses the main properties that a Ritt characteristic set of $\langle\pset{P}\rangle$ has.

\begin{proposition}\label{charpro}
Let $\pset{C}=[C_1,\ldots, C_r]$ be the W-characteristic set of $\langle\pset{P}\rangle$ with $I_{i}=\ini(C_{i})$ for $1\leq i\leq r$. Then:
\begin{enumerate}
\item[{\rm (a)}] $\prem(P, \pset{C})\equiv 0$ for all $P\in \langle\pset{P}\rangle$;
\item[{\rm (b)}] $\langle\pset{C}\rangle\subset \langle\pset{P}\rangle\subset\sat(\pset{C})$;
\item[{\rm (c)}] $\zero(\pset{C})\setminus\zero(\{I_1\cdots I_r\})\subset\zero(\pset{P})
    \subset\zero(\pset{C})$.
\end{enumerate}
\end{proposition}
\emph{Proof.} (a) Let $p_i=\cls(C_i)$ for $1\leq i\leq r$. Then $\pset{G}=\pset{G}^{\langle p_1\rangle}\cup\cdots\cup\pset{G}^{\langle p_r\rangle}$ and for all $C\in\pset{G}^{\langle p_i\rangle}\setminus\{C_i\}$, $C>_{{\rm plex}} C_i$ and thus $\deg(\lpp(C), x_{p_i})=\deg(C, x_{p_i})\geq \deg(C_i, x_{p_i})$.
Hence, for any $P\in \langle\pset{P}\rangle$, $R=\prem(P, \pset{C})$ is B-reduced with respect to $\pset{G}$. Therefore, $R\equiv 0$ and (a) is proved.

(b) Note that $\pset{C}\subset\pset{G}\subset \langle\pset{P}\rangle$, so $\langle\pset{C}\rangle\subset \langle\pset{P}\rangle$. For any $P\in \langle\pset{P}\rangle$,
by (a) and the pseudo-remainder formula for $\prem(P, \pset{C})=0$ there exist nonnegative  integers $q_1,\ldots,q_r$ such that
$I_1^{q_1}\cdots I_r^{q_r}P\in \langle\pset{C}\rangle$. It follows from the definition of saturated ideals that $P\in\sat(\pset{C})$. Therefore, (b) is proved.

(c) $\zero(\pset{P}) \subset\zero(\pset{C})$ follows from the first $\subset$ relation in (b).
Consider any zero $\p{a}\in \zero(\pset{C})\setminus\zero(\{I_1\cdots I_r\})$ and let $P$ be any polynomial in $\pset{P}$. Then $C_{i}(\p{a})=0$ and $I_i(\p{a})\neq 0$ for $1\leq i\leq r$. Plunging $\p{a}$ into the pseudo-remainder formula for $\prem(P, \pset{C})=0$, we see that $P(\p{a})=0$.
Therefore, $\p{a}\in\zero(\pset{P})$ and (c) is proved. $\square$

\begin{proposition}\label{elimpro}
Let $\pset{C}$ be the W-characteristic set of $\langle\pset{P}\rangle$. Then for every $i$ $(0\leq i\leq n)$, $\pset{C}^{\langle 0,\ldots,i\rangle}$ is the W-characteristic set of the $(n-i)$th elimination ideal $\langle\pset{P}\rangle^{\langle 0,\ldots,i\rangle}$ of $\langle\pset{P}\rangle$.
\end{proposition}
\emph{Proof.} By the elimination theorem of Buchberger-Gr\"obner bases (see, e.g., \cite[Ch.\,3, Theorem~2]{clo96}), $\pset{G}^{\langle 0,\ldots,i\rangle}$ is the reduced plex Buchberger-Gr\"obner basis of $\langle\pset{P}\rangle^{\langle 0,\ldots,i\rangle}$. Hence the W-characteristic set of $\langle\pset{P}\rangle^{\langle 0,\ldots,i\rangle}$ is identical to $\pset{C}^{\langle 0,\ldots,i\rangle}$. $\square$

\begin{theorem}\label{mainth0}
Let $\pset{C}$ be the W-characteristic set of $\langle\pset{P}\rangle$.
If $\pset{C}$ is an ascending set, then $\pset{C}$ is a Ritt characteristic set of $\langle\pset{P}\rangle$.
\end{theorem}
\emph{Proof.} By Proposition~\ref{charpro} (a), $\prem(P, \pset{C})\equiv 0$ for all $P\in \langle\pset{P}\rangle$. The theorem follows from Lemma~\ref{rittlem}.
$\square$

\begin{theorem}\label{rittcsthm}
Let $\pset{C}=[C_1,\ldots,C_r]$ be the W-characteristic set of $\langle\pset{P}\rangle$. If $\pset{C}$ is regular,
then \[\pset{C}^*=[C_1, \prem(C_2, [C_1]), \ldots, \prem(C_r, [C_1,\ldots,C_{r-1}])]\]
is a Ritt characteristic set of $\langle\pset{P}\rangle$, where $\pset{C}^*$ is also regular.
\end{theorem}
\emph{Proof.} Let $\pset{C}_i=[C_1,\ldots,C_i]$, $x_{p_i}=\lv(C_i)$, $I_i=\ini(C_i)$, and $\pset{C}_i^*=[C_1^*,\ldots,C_i^*]$ for $1\leq i\leq r$, where $C_1^*=C_1$ and $C_i^*=\prem(C_i, \pset{C}_{i-1})$ for $2\leq i\leq r$. Note first that $C_i^*\in\langle \pset{C}\rangle$, so $\langle \pset{C}^*\rangle\subset \langle \pset{C}\rangle$. Let $J$ be any power product of $I_1^*=\lc(C_1^*,x_{p_1}),\ldots,I_r^*=\lc(C_r^*,x_{p_r})$. Observe from the pseudo-remainder formula for $C_i^*=\prem(C_i, \pset{C}_{i-1})$ that for each $I_i^*$, there exist a power product $J_i$ of $I_1,\ldots,I_i$ and polynomials $Q_{i,1},\ldots,Q_{i,i-1}$ in $\fk[x_1,\ldots,x_{n}]$ such that \[J_i-(Q_{i,1}C_1+\cdots+Q_{i,i-1}C_{i-1})=I_i^*.\] Multiplying the two sides of such equalities up to certain powers, one sees that $J$ is equal to a power product $H$ of $J_1,\ldots,J_r$ plus a linear combination $F=Q_1C_1+\cdots+Q_{r-1}C_{r-1}$ for some polynomials $Q_i\in\fk[x_1,\ldots,x_n]$, i.e., $J=H+F$, so $H-J\in\langle\pset{C}\rangle$.
As $H$ is also a power product of $I_1,\ldots,I_r$ and $\pset{C}$ is regular, $N=\res(H,\pset{C})\not\equiv 0$ and $N$ does not involve $x_{p_1}, \ldots, x_{p_r}$. According to \eqref{resform} there exists a polynomial $A$ in $\fk[x_1,\ldots,x_n]$ such that $N-AH\in\langle\pset{C}\rangle$. This implies that $N-AJ\in\langle\pset{C}\rangle$.
Therefore, $J\not\equiv 0$; for otherwise, $N=\prem(N, \pset{C})=\prem(N-AJ,\pset{C})\equiv 0$ (by Proposition~\ref{satpro}) leads to contradiction. In particular, we have $I_i^*\not\equiv 0$ for all $i$. Since $\deg(C_i^*,x_{p_i})$ cannot be greater than $\deg(C_i,x_{p_i})$, they must be equal. This shows that $\lv(C_i^*)=x_{p_i}$ and $\ini(C_i^*)=I_i^*$ for $1\leq i\leq r$. Consequently, $C_i^*$ is R-reduced with respect to $\pset{C}_{i-1}^*$, as it is so with respect to $\pset{C}_{i-1}$, for $2\leq i\leq r$; thereby $\pset{C}^*$ is an ascending set.

For any polynomial $P\in\sat(\pset{C}^*)$, there exists a power product $J$ of $I_1^*,\ldots,I_r^*$ such that
 $JP\in\langle\pset{C}^*\rangle\subset\langle\pset{C}\rangle$. According to the above reasoning, we have $N-AJ\in\langle\pset{C}\rangle$. It follows that $NP\in\langle\pset{C}\rangle$.
On the other hand, let $R^*=\prem(P, \pset{C}^*)$. Then there exists a power product $K$ of $I_1^*,\ldots,I_r^*$ such that $KP-R^*\in\langle\pset{C}^*\rangle\subset\langle\pset{C}\rangle$. It follows that $NR^*\in\langle\pset{C}\rangle$.
Obviously, $\deg(R^*, x_{p_i})<\deg(C_i^*,x_{p_i})=\deg(C_i,x_{p_i})$ for all $i$. Hence $NR^*=\prem(NR^*, \pset{C})\equiv 0$. Therefore, $R^*\equiv 0$ and thus $\pset{C}^*$ is regular by Proposition~\ref{satpro}.

 Next we want to show that $\langle\pset{C}\rangle\subset\sat(\pset{C}^*)$. For this purpose, assume by induction that $C_{i-1}\in \sat(\pset{C}_{i-1}^*)$; the case for $C_1$ is trivial.
 Then there exists a power product $K_i$ of $I_1, \ldots, I_{i-1}$ such that $K_iC_i-C_i^*\in\langle \pset{C}_{i-1}\rangle\subset\sat(\pset{C}_{i-1}^*)$. Since $\pset{C}$ is regular, $R_i=\res(K_i,\pset{C}_{i-1})\not\equiv 0$ and $R_i$ does not involve $x_{p_1}, \ldots, x_{p_r}$ for each $i$. According to \eqref{resform} there exists a polynomial $A_i$ in $\fk[x_1,\ldots,x_n]$ such that $R_i-A_iK_i\in\langle\pset{C}_{i-1}\rangle\subset\sat(\pset{C}_{i-1}^*)$.
This implies that $R_iC_i\in\sat(\pset{C}_{i}^*)$. It follows from Proposition~\ref{satpro} that \[R_i\prem(C_i,\pset{C}_{i}^*)=\prem(R_iC_i,\pset{C}_{i}^*)\equiv 0.\]
 Hence $\prem(C_i,\pset{C}_{i}^*)\equiv 0$ and $C_{i}\in \sat(\pset{C}_{i}^*)$. Therefore, $\langle\pset{C}\rangle\subset\sat(\pset{C}^*)$.

 Now consider any polynomial $P\in\langle\pset{P}\rangle$. Clearly, $R=\prem(P, \pset{C})$ is B-reduced
with respect to the Buchberger-Gr\"obner basis $\pset{G}$ of $\pset{P}$. This implies that $R\equiv 0$ and
$P\in\sat(\pset{C})$. Hence there exists a power product $L$ of $I_1, \ldots, I_{r}$ such that $LP\in\langle \pset{C}\rangle$. Since $M=\res(L, \pset{C})\not\equiv 0$, there exists a polynomial $B$ in $\fk[x_1,\ldots,x_n]$ such that $M-BL\in\langle\pset{C}\rangle$. It follows that $MP\in\langle \pset{C}\rangle\subset\sat(\pset{C}^*)$.  As $M$ does not involve $x_{p_1},\ldots,x_{p_r}$, $M\prem(P, \pset{C}^*)=\prem(MP, \pset{C}^*)\equiv 0$ according to Proposition~\ref{satpro}; therefore
$\prem(P, \pset{C}^*)\equiv 0$. By Lemma~\ref{rittlem},
$\pset{C}^*$ is a Ritt characteristic set of $\langle\pset{P}\rangle$.
$\square$

\begin{corollary}\label{maincor0}
Let $\pset{C}$ be the W-characteristic set of $\langle\pset{P}\rangle$.
For any $0\leq i\leq n$, if $\pset{C}^{\langle 0,\ldots,i\rangle}=[C_1,\ldots,C_k]$ is regular, then
\[[C_1, \prem(C_2, [C_1]), \ldots, \prem(C_k, [C_1,\ldots,C_{k-1}])]\]
is a Ritt characteristic set of the elimination ideal $\langle\pset{P}\rangle^{\langle 0,\ldots,i\rangle}$.
\end{corollary}
\emph{Proof.} It follows from Proposition~\ref{elimpro} and Theorem~\ref{rittcsthm}.
$\square$

\subsection{Structure of Buchberger-Gr\"obner bases}

 To explore the structural properties of the W-characteristic set $\pset{C}$ of $\langle\pset{P}\rangle$, we assume from now on that $x_1,\ldots,x_n$ are properly ordered such that the leading variables $x_{p_i}$ of the polynomials in $\pset{C}$ are greater than all the other free variables, called \emph{parameters}. Let $y_1=x_{p_1}, \ldots, y_r=x_{p_r}$ and $u_1,\ldots,u_m$ be all the parameters, where $m+r=n$. Then the assumed variable order is
 $u_1<_{{\rm plex}}\cdots<_{{\rm plex}}u_m<_{{\rm plex}}y_1<_{{\rm plex}}\cdots<_{{\rm plex}}y_r$.

\begin{theorem}\label{mainth1}
Let $[C_1,\ldots, C_r]$ be the W-characteristic set of $\langle\pset{P}\rangle$. For any $1\leq k<r$, if $\pset{C}_k=[C_1,\ldots,C_k]$ is normal and $I_{k+1}=\ini(C_{k+1})$, with $\lv(I_{k+1})=y_l$, is not R-reduced with respect to $C_l$, then \[\prem(I_{k+1}, \pset{C}_k)\equiv 0\quad\mbox{and}\quad \prem(C_{k+1}, \pset{C}_k)\equiv 0.\]
\end{theorem}
\emph{Proof.} Let $y_i=\lv(C_i)$ and $I_i=\ini(C_i)$ for $1\leq i\leq r$.
Suppose that $I_{k+1}$ is not R-reduced with respect to
$C_{l}$ and let $R=\prem(I_{k+1}, \pset{C}_k)$. Then according to \eqref{rmdr:TS} there exists a power product $T$ of $I_1,\ldots, I_k$ such that $TI_{k+1}-R=D\in\langle \pset{C}_k\rangle$. As $\lv(I_{k+1})=y_l$, $I_{k+1}$ does not involve $y_{{l+1}},\ldots, y_{k}$. This implies that
$\deg(R, y_{l})<\deg(I_{k+1},y_{l})$ and $\deg(R, y_{i})=\deg(I_{k+1},y_{i})=0$ for $l+1\leq i\leq k$. It follows that $R<_{{\rm plex}}I_{k+1}$ and thus $Ry_{{k+1}}^d<_{{\rm plex}}I_{k+1}y_{{k+1}}^d$, where $d=\deg(C_{k+1},y_{{k+1}})$. Multiplying the two sides of
$C_{k+1}=I_{k+1}y_{{k+1}}^d+H_{k+1}$ by $T$, we have
\[TC_{k+1}=TI_{k+1}y_{{k+1}}^d+TH_{k+1}=Ry_{{k+1}}^d+Dy_{{k+1}}^d+TH_{k+1}
\in\langle\pset{P}\rangle.\]
Since $\deg(H_{k+1},y_{{k+1}})<d$ and $T$ does not involve $y_{{k+1}}$,
$TH_{k+1}$ is B-reduced with respect to $\pset{G}^{\langle m+k+1\rangle}$. Note that $D\in\langle\pset{C}_k\rangle\subset\langle\pset{G}\rangle$, so $\nf(D y_{{k+1}}^d, \pset{G})\equiv 0$.
If $R\not\equiv 0$, then \[0\neq R y_{{k+1}}^d+\nf(TH_{k+1}, \pset{G})=R y_{{k+1}}^d+\nf(D y_{{k+1}}^d+T H_{k+1}, \pset{G})\in\langle\pset{P}\rangle\]
is B-reduced with respect to $\pset{G}^{\langle 0, \ldots, {m+k+1}\rangle}$, which leads to contradiction. Therefore, $\prem(I_{k+1}$, $\pset{C}_{k})=R\equiv 0$ and $\nf(TH_{k+1}, \pset{G})\equiv 0$.

Moreover, $\nf(TH_{k+1}, \pset{G})\equiv 0$ implies that $\prem(TH_{k+1},\pset{C}_k)\equiv 0$.
 Since $\pset{C}_k$ is normal, $T$ does not involve $y_{1},\ldots,y_{{k}}$ and $\prem(TI_{k+1} y_{{k+1}}^d, \pset{C}_k)=Ty_{{k+1}}^d\prem(I_{k+1}, \pset{C}_k)\equiv 0$. It follows from Corollary~\ref{satcor} that \[T\prem(C_{k+1}, \pset{C}_{k})=\prem(TC_{k+1}, \pset{C}_{k})=\prem(TI_{k+1} y_{{k+1}}^d+TH_{k+1}, \pset{C}_{k})\equiv 0.\]
 Therefore, $\prem(C_{k+1}, \pset{C}_{k})\equiv 0$ and the proof is complete.
$\square$

\begin{lemma}\label{mainlem}
Let $[C_1,\ldots, C_r]$ be
the W-characteristic set of $\langle\pset{P}\rangle$ with $\pset{C}_{i}=[C_1,\ldots,C_i]$ for $1\leq i\leq r$ and let $k$ be the biggest integer such that
$\pset{C}_{k}$ is normal. Assume that $k<r$ and let $I_{k+1}=\ini(C_{k+1})$, $y_l=\lv(I_{k+1})$, and $Q=\pquo(C_{l}, I_{k+1}, y_{{l}})$. Then:
\begin{enumerate}
\item[{\rm (a)}] $\pset{C}_{k+1}$ is not regular;

\item[{\rm (b)}] if $I_{k+1}$ is R-reduced with respect to $C_{l}$, then \[\prem(C_{l}, [C_1, \ldots, C_{l-1},I_{k+1}])\equiv 0\quad
    \mbox{and}\quad \prem(QC_{k+1}, \pset{C}_{k})\equiv 0.\]
\end{enumerate}
\end{lemma}
\emph{Proof.}
(a) Let $y_{i}=\lv(C_i)$ for $1\leq i\leq r$ and $R=\res(I_{k+1}, \pset{C}_{l})$.
Then according to \eqref{resform} there exist polynomials $A, B_1, \ldots, B_{l}\in\fk[u_1,\ldots,u_m, y_1,\ldots,y_{l}]$ such that
\begin{equation}\label{resfor}
R=AI_{k+1}+B_1C_1+\cdots+B_{l}C_{l},
\end{equation}
where $R$ does not involve $y_{1},\ldots,y_{l}$.
Write $C_{k+1}=I_{k+1} y_{{k+1}}^d+H_{k+1}$, where $d=\deg(C_{k+1},y_{{k+1}})$. Multiplying the two sides of this equality
by $A$ and using \eqref{resfor}, one obtains
\[R y_{{k+1}}^d+AH_{k+1}=AC_{k+1}+(B_1C_1+\cdots+B_{l}C_{l}) y_{{k+1}}^d\in\langle\pset{P}\rangle.\]
Suppose that $\pset{C}_{k+1}$ is regular. Then $R\not\equiv 0$. Note that $R$ does not involve
$y_{1},\ldots,y_{r}$, so $R$ is R-reduced with respect to
$\pset{C}_k$ and B-reduced with respect to $\pset{G}^{\langle 0,\ldots,m+k\rangle}$; thereby $R <_{{\rm plex}} I_{k+1}$ and $R y_{{k+1}}^d <_{{\rm plex}} I_{k+1} y_{{k+1}}^d$. Hence
$R y_{{k+1}}^d$ is B-reduced with respect to $\pset{G}^{\langle 0,\ldots,m+k+1\rangle}$.
It follows that
\[0\neq R y_{{k+1}}^d+\nform(AH_{k+1},\pset{G})\in\langle\pset{P}\rangle\] is B-reduced with respect to $\pset{G}$. This leads to contradiction. Therefore, $R\equiv 0$ and $\pset{C}_{k+1}$ is not regular.

(b) Suppose that $I_{k+1}$ is R-reduced with respect to $C_{l}$ and let \[M=\prem(\prem(C_{l},I_{k+1},y_{{l}}), \pset{C}_{l-1}),\]
 $I=\ini(I_{k+1})$, and $I_i=\ini(C_i)$ for $1\leq i\leq k$. Then there exist an integer $q\geq 0$ and a power product $U$ of $I_1, \ldots, I_{l-1}$ such that
\begin{equation}\label{uqceq}
U(I^{q}C_{l} - QI_{k+1})-M=E\in\langle\pset{C}_{l-1}\rangle.
\end{equation}
Recall that $C_{k+1}=I_{k+1} y_{{k+1}}^d+H_{k+1}$. It follows that
\[UQC_{k+1}=(I^{q}UC_{l}-M) y_{{k+1}}^d-E y_{{k+1}}^d+UQH_{k+1}\in\langle\pset{P}\rangle.\]
As $I_{k+1}$ is R-reduced with respect to $C_{l}$,  $\deg(I_{k+1},y_{{l}})<\deg(C_{l},y_{{l}})$.
Note that $M$ does not involve $y_{{l+1}}, \ldots, y_{k}$, so $M$ is R-reduced with respect to $\pset{C}_k$ and B-reduced with respect to $\pset{G}^{\langle 0,\ldots,m+k\rangle}$. Moreover, $M$ is R-reduced with respect to $I_{k+1}$ and thus $M <_{{\rm plex}} I_{k+1}$ and $M y_{{k+1}}^d <_{{\rm plex}} I_{k+1} y_{{k+1}}^d$. Hence
$M y_{{k+1}}^d$ is B-reduced with respect to $\pset{G}^{\langle 0,\ldots,m+k+1\rangle}$.
If $M\not\equiv 0$, then \[\begin{array}{ll}\smallskip
0\!\!\!& \neq -M y_{{k+1}}^d+\nf(UQH_{k+1}, \pset{G})\\
&=-M y_{{k+1}}^d+\nf(I^{q}UC_{l} y_{{k+1}}^d-E y_{{k+1}}^d+UQH_{k+1}, \pset{G})\in\langle\pset{P}\rangle\end{array}\] is B-reduced with respect to $\pset{G}$, which leads to contradiction. Therefore, $M$ must be identically equal to $0$.

Since $\pset{C}_k$ is normal, $I^qUC_l\in\langle\pset{C}_l\rangle$, and $E\in\langle\pset{C}_{l-1}\rangle$, we have $I^qUC_l-E\in\langle\pset{C}_l\rangle$ and
 \[\prem(UQI_{k+1} ,\pset{C}_l)=\prem(I^{q}UC_{l} -E, \pset{C}_l)\equiv 0 .\]
It follows that $\prem(UQI_{k+1}y_{k+1}^d,\pset{C}_k)\equiv 0$ and $UQH_{k+1}\in\langle\pset{P}\rangle$. As $\deg(H_{k+1}, y_{k+1})<d$, $UQH_{k+1}$ is R-reduced with respect to $C_{k+1}$; thereby \[\prem(UQH_{k+1}, \pset{C}_{k})=\prem(UQH_{k+1}, \pset{C}_{k+1})\equiv 0.\]  Hence $\prem(UQC_{k+1}, \pset{C}_{k})\equiv 0$ by Corollary~\ref{satcor}.
 As $U$ does not involve $y_{1},\ldots,y_{{k}}$, \[\prem(UQC_{k+1}, \pset{C}_{k})=U\prem(QC_{k+1}, \pset{C}_{k}).\]
 Therefore, $\prem(QC_{k+1}, \pset{C}_{k})\equiv 0$. The proof is complete.
$\square$

\begin{lemma}\label{mainlem1}
Let $[C_1,\ldots, C_r]$ be
the W-characteristic set of $\langle\pset{P}\rangle$ and $k$ be the biggest integer such that
$[C_1,\ldots,C_k]$ is normal. Assume that $k<r$ and let $I_{k+1}=\ini(C_{k+1})$ and $y_l=\lv(I_{k+1})$. If $I_{k+1}$ is R-reduced with respect to $C_{l}$, then either $\res(\ini(I_{k+1}), \pset{C}_{l-1})\equiv 0$, or
 \[\prem(C_{k+1}, [C_1, \ldots, C_{l-1},I_{k+1},C_{l+1},\ldots,C_k])\equiv 0.\]
\end{lemma}
\emph{Proof.} Suppose that $I_{k+1}$ is R-reduced with respect to $C_{l}$,
let $y_i=\lv(C_i)$, $I_i=\ini(C_i)$, and $\pset{C}_{i}=[C_1,\ldots,C_i]$ for $1\leq i\leq r$, and let
\[\tilde{\pset{C}}=[C_1, \ldots, C_{l-1},I_{k+1},C_{l+1},\ldots,C_k],\quad d=\deg(C_{k+1}, y_{k+1}),\quad
H_{k+1}=C_{k+1}-I_{k+1} y_{{k+1}}^d,\] and
$R=\prem(H_{k+1}, \tilde{\pset{C}})$. Then there exist nonnegative integers $s$, $s_{1},\ldots,s_{l-1}$, $s_{l+1},\ldots,s_k$ and polynomials
$B, B_{1},\ldots,B_{l-1},B_{l+1},\ldots,B_k\in\fk[u_1,\ldots,u_m,y_1,\ldots,y_{k}]$
such that
\begin{equation}\label{qreq}
I^sI_{1}^{s_{1}}\cdots I_{l-1}^{s_{l-1}} I_{l+1}^{s_{l+1}}\cdots I_k^{s_k}QH_{k+1}
= BQI_{k+1} - QB_{l}C_{l} +\sum_{i=1}^kQB_{i}C_{i}+QR,\end{equation}
where $Q=\pquo(C_{l}, I_{k+1}, y_{{l}})$ and $I=\ini(I_{k+1})$.
The first conclusion of Lemma~\ref{mainlem} (b) implies that $L=\prem(C_{l}, I_{k+1}, y_{{l}})\in\sat(\pset{C}_{l-1})$, i.e., $I^qC_{l}-QI_{k+1}=L\in\sat(\pset{C}_{l-1})$ for some integer $q\geq 0$. Recall that $\pset{C}_k$ is normal. Hence $QI_{k+1}\in\sat(\pset{C}_{l})$ and $\prem(QI_{k+1}, \pset{C}_{l})\equiv 0$.
This, together with the second conclusion of Lemma~\ref{mainlem} (b), implies that $\prem(QH_{k+1}, \pset{C}_{k})=\prem(QC_{k+1}-QI_{k+1} y_{{k+1}}^d, \pset{C}_{k})\equiv 0$, so that
$QH_{k+1}\in\sat(\pset{C}_{k})$. It follows from \eqref{qreq} that $QR\in\sat(\pset{C}_{k})$.

As $I_{k+1}$ is R-reduced with respect to $C_{l}$, $\deg(Q,y_{{l}})>0$. From the pseudo-remainder formula $I^qC_{l}=QI_{k+1}+L$, one sees that $\ini(Q)=I^{q-1}I_{l}$ for some $q\geq 0$.
Note that $I_{l}$ does not involve $y_{1},\ldots, y_{{l-1}}$ and assume that $\res(I, \pset{C}_{l-1})\not\equiv 0$. Then $M=\res(I^{q-1}I_{l}, \pset{C}_{l-1})\not\equiv 0$.
Hence there exists a polynomial $S\in\fk[u_1,\ldots,u_m,y_1,\ldots,y_{{l-1}}]$
such that $M-S\ini(Q)=A\in\langle\pset{C}_{l-1}\rangle$. Write $Q=\ini(Q)y_{l}^{\delta}+\bar{Q}$ and let $Z=My_{l}^{\delta}+S\bar{Q}$, where $\delta=\deg(Q,y_{l})$. Then $SQ=Z-Ay_{l}^{\delta}$ and thus $SQR=ZR-ARy_{l}^{\delta}$. It follows that $ZR\in\sat(\pset{C}_{k})$. Since $\pset{C}_{k}$ is normal,
$\prem(ZR, \pset{C}_{k})\equiv 0$ according to Proposition~\ref{satpro}.
It is easy to see that
\[\begin{array}{ll}\smallskip
\deg(ZR,y_{l})\!\!\!&=\deg(Z,y_{l})+\deg(R,y_{l})=\delta+\deg(R,y_{l})\\
&=\deg(C_l,y_{l})-\deg(I_{k+1},y_{l})+\deg(R,y_{l})<\deg(C_l,y_{l}).
\end{array}\]
Hence \[\prem(MRy_{l}^{\delta}+S\bar{Q}R, \pset{C}_{l-1})=\prem(ZR, \pset{C}_{l-1})=\prem(ZR, \pset{C}_{l})=\prem(ZR, \pset{C}_{k})\equiv 0\] (for $Z$ does not involve $y_{{l+1}},\ldots,y_{k}$). This implies that $MR=\prem(MR, \pset{C}_{l-1})\equiv 0$ according to Lemma~\ref{satlem}.
Therefore, $R\equiv 0$.

As $\pset{C}_k$ is normal and $\res(I, \pset{C}_{l-1})\not\equiv 0$, $\tilde{\pset{C}}$ is regular.
Obviously, $\prem(I_{k+1} y_{{k+1}}^d,\tilde{\pset{C}})=y_{{k+1}}^d\prem(I_{k+1},\tilde{\pset{C}})\equiv 0$. Hence \[\prem(C_{k+1}, \tilde{\pset{C}})=\prem(I_{k+1} y_{{k+1}}^d+H_{k+1}, \tilde{\pset{C}})\equiv 0\] by Corollary~\ref{satcor}. The proof is complete.
$\square$

\begin{theorem}\label{mainth}
Let $\pset{C}=[C_1,\ldots,C_r]$ be the W-characteristic set of $\langle\pset{P}\rangle$ and
$\pset{C}_{i}=[C_1,\ldots,C_i]$ for $1\leq i\leq r$.
If $\pset{C}$ is abnormal, then there exists an integer $k$ $(1\leq k<r)$
such that
\begin{enumerate}
\item[{\rm (a)}] $\pset{C}_{k}$ is normal and thus regular;
\item[{\rm (b)}] $\pset{C}_{k+1}$ is not regular;
\item[{\rm (c)}] if $I_{k+1}=\ini(C_{k+1})$ is not R-reduced with respect to $C_{l}$, then
\[\prem(I_{k+1}, \pset{C}_{l})\equiv 0\quad\mbox{and}\quad \prem(C_{k+1}, \pset{C}_{k})\equiv 0;\]

 \item[{\rm (d)}] if $I_{k+1}$ is R-reduced with respect to $C_{l}$, then $\prem(C_{l}, [C_1, \ldots, C_{l-1},I_{k+1}])\equiv 0$ and either $\res(\ini(I_{k+1}), \pset{C}_{l-1})\equiv 0$, or
 \[\prem(C_{k+1}, [C_1, \ldots, C_{l-1},I_{k+1},C_{l+1},\ldots,C_k])\equiv 0,\]
\end{enumerate}
where $y_{{l}}=\lv(I_{k+1})$.
\end{theorem}
\emph{Proof.} Suppose that $\pset{C}$ is abnormal and let $k$ ($1\leq k<r$) be the biggest integer such that $\pset{C}_{k}$ is normal. Then we have (a).
By Lemma~\ref{mainlem} (a), $\pset{C}_{k+1}$
is not regular; thus (b) is proved. The identity $\prem(C_{k+1}, \pset{C}_{k})\equiv 0$ in (c) has been proved as the second conclusion of Theorem~\ref{mainth1}. Recall $\prem(I_{k+1}, \pset{C}_k)\equiv 0$, the first conclusion of Theorem~\ref{mainth1}, where $I_{k+1}=\ini(C_{k+1})$. As $y_{{l+1}},\ldots,y_{k}$ do not appear in $I_{k+1}$, one sees that $\prem(I_{k+1}, \pset{C}_k)=\prem(I_{k+1}, \pset{C}_l)$. Hence $\prem(I_{k+1}, \pset{C}_{l})\equiv 0$ and (c) is proved.

The identity $\prem(C_{l}, [C_1, \ldots, C_{l-1},I_{k+1}])\equiv 0$ in (d) has been proved as the first conclusion of Lemma~\ref{mainlem} (b), and so has the second conclusion of (d) proved as Lemma~\ref{mainlem1}.
$\square$

\medskip
The irregularity index ${m+k+1}$ in Theorem~\ref{mainth} is clearly characteristic. The W-characteristic set of the $(n-m-k-1)$th elimination ideal $\langle\pset{P}\rangle^{\langle 0,\ldots,{m+k+1}\rangle}$ is irregular, while that of the $(n-m-k)$th elimination ideal $\langle\pset{P}\rangle^{\langle 0,\ldots,{m+k}\rangle}$
is not only regular but also normal. As shown in Corollary~\ref{maincor0}, from the normal W-characteristic set of the elimination ideal a Ritt characteristic set of the ideal can be computed rather easily by means of pseudo-division. When the W-characteristic set of the ideal $\langle\pset{P}\rangle$ is itself normal, ${m+k+1}$ may be defined to be $n+1$ (or any other integer greater than $n$) and the Buchberger-Gr\"obner basis $\pset{G}$ of $\pset{P}$ may be said to be \emph{regular}.

The pseudo-divisibility relations in Theorem~\ref{mainth} (c) and (d) expose the intrinsic structure of the polynomials $C_{l}$ and $C_{k+1}$ and the irregularity of $C_{k+1}$ modulo the saturated ideal of the normal triangular set $\pset{C}_k$. More relations of this kind would help us gain more insights into the structure of the polynomials in $\pset{G}$.

\subsection{Examples}

The following examples serve to illustrate various behaviors of the W-characteristic sets of polynomial ideals.

\begin{example}\label{exa}\em
(a) Let $\pset{P}=\{x_1x_2-1, x_3-x_2\}$. Then the W-characteristic set of $\langle\pset{P}\rangle$ is $\pset{C}=[x_1x_2-1, x_3-x_2]$. $\pset{C}$ is normal, but it is not a Ritt characteristic set of $\langle\pset{P}\rangle$.
Construct $\pset{C}^*=[x_1x_2-1, \prem(x_3-x_2,[x_1x_2-1])]=[x_1x_2-1, x_1x_3-1]$. Then
$\pset{C}^*$ is a Ritt characteristic set of $\langle\pset{P}\rangle$ and $\pset{C}^*\sim\pset{C}$, but $x_1x_3-1>_{{\rm plex}} x_3-x_2$.

\medskip
(b) Let $\pset{P}=\{x_1^2, (x_2+x_1)x_3+x_1\}$. Then the W-characteristic set of $\langle\pset{P}\rangle$ is $\pset{C}=[x_1^2, (x_2+x_1)x_3+x_1]$: it is a regular ascending set and thus is a Ritt characteristic set of $\langle\pset{P}\rangle$. The regular set $\pset{C}$ is not normal because the parameter $x_2$ is ordered greater than the leading variable $x_1$. This example explains why the assumption on the variable order is necessary for the W-characteristic set of $\langle\pset{P}\rangle$ to be normal or exhibit irregularity structure.

\medskip
(c) Let $\pset{P}=\{x_1x_2, x_2x_3, x_3x_4\}$. Then the W-characteristic set of $\langle\pset{P}\rangle$ is $\pset{C}=[C_1,C_2,C_3]=[x_1x_2, x_2x_3, x_3x_4]$.
$\pset{C}$ is abnormal and irregular and it is not a Ritt characteristic set of $\langle\pset{P}\rangle$. One may see that $\prem(C_2, C_1, x_2)=\prem(C_3, C_2, x_3)\equiv 0$ and $[x_1x_2, x_3x_4]$ is a Ritt characteristic set of $\langle\pset{P}\rangle$. However,
application of the construction for $C_i^*$ in Theorem~\ref{rittcsthm} to the irregular W-characteristic set $\pset{C}$ here does not lead to the Ritt characteristic set $[x_1x_2, x_3x_4]$ of $\langle \pset{P}\rangle$.

\medskip
(d) Let $C_3=x_1(x_2+x_1)x_3-x_1^3$ and $\pset{P}=\{x_1^4, x_2^4, C_3\}$. Then $\pset{G}=\{x_1^4, x_1^3x_2^3, x_2^4, C_3\}$. Thus the W-characteristic set of $\langle\pset{P}\rangle$ is $\pset{C}=[x_1^4, x_1^3x_2^3, C_3]$, which is abnormal and irregular. By Theorem~\ref{mainth0}, $\pset{C}$ is a Ritt characteristic set of $\langle\pset{P}\rangle$.

Now let $C_1=x_1^3$, $C_2=x_2^3$, $\bar{C}_3=x_1x_2x_3-x_1^2x_2$ and $\bar{\pset{P}}=\{C_1, C_2, \bar{C}_3\}$. Then the W-characteristic set $\bar{\pset{C}}=[C_1, C_2, \bar{C}_3]$ of $\langle\bar{\pset{P}}\rangle$ is a Ritt characteristic set of $\langle\bar{\pset{P}}\rangle$ by Theorem~\ref{mainth0}. One can see that $Q=\pquo(C_2, \ini(\bar{C}_3), x_2)=x_1^2x_2^2$ and $\res(\ini(Q), [C_1])\equiv 0$. This example shows that the case in which $\res(I, \pset{C}_{l-1})\equiv 0$ in Theorem~\ref{mainth} (d) does occur.

\medskip
(e) Let $C_1=x_1x_2$, $C_2=x_3x_4-x_2^2$, $C_3=x_2x_5+x_4^2$ and $\pset{P}=\{C_1, C_2, C_3\}$. Then $\pset{G}=\{C_1, C_2, x_1x_4^2, C_3\}$ and the W-characteristic set of $\langle\pset{P}\rangle$ is $\pset{C}=[C_1, C_2, C_3]$. One can verify that $\prem(C_3, [C_1, C_2])\equiv 0$ and $\prem(x_1x_5-x_1, \pset{C})\equiv 0$, but $x_1x_5-x_1\not\in\langle\pset{P}\rangle$.

Let $\bar{C}_2=x_2x_4-x_2^2$ instead of $C_2$ and $\bar{\pset{P}}=\{C_1, \bar{C}_2, C_3\}$. Then $\bar{\pset{G}}=\{C_1, \bar{C}_2, x_1x_4^2, x_4^3-x_2^3, C_3\}$ and the W-characteristic set of $\langle\bar{\pset{P}}\rangle$ is $\bar{\pset{C}}=[C_1, \bar{C}_2, C_3]$. Now $\prem(C_3, [C_1, \bar{C}_2])=\prem(\bar{C}_2,[C_1])\equiv 0$. However,
$x_1x_4^2\in\langle\bar{\pset{P}}\rangle$, but $\prem(x_1x_4^2, [C_1])\not\equiv 0$. Therefore, $[C_1]$ is not a Ritt characteristic set of $\langle\bar{\pset{P}}\rangle$, and we suspect that $[C_1, x_1x_4^2]$ is.
 This example shows that, in the abnormal case, the minimal ascending set contained in the W-characteristic set $\pset{C}$ of an ideal $\ideal{I}$ is not necessarily a Ritt characteristic set of the ideal and the Buchberger-Gr\"obner basis $\pset{G}$ of $\ideal{I}$ may contain other ascending sets of lower rank. A natural question that remains to be answered is how to construct a Ritt characteristic set of $\ideal{I}$ from $\pset{G}$ when $\pset{C}$ is neither an ascending set nor a regular set.

\medskip
For Example~\ref{exa} (c), as well as other examples we have studied, it appears that the Ritt characteristic set of an ideal $\ideal{I}$ does not characterize the ideal well enough in the abnormal case.
This is caused essentially by the irregularity of the ideal. The W-characteristic set of $\ideal{I}$, which provides sufficient information about $\ideal{I}$, may serve as an alternative to the Ritt characteristic set.
\end{example}

\section{Some Remarks}

We point out two directions of research in which triangular sets and Buchberger-Gr\"obner bases have
been explored reciprocally. The first is concerned with algorithmic decomposition of
polynomial or differential polynomial sets (or systems) into triangular or differential
triangular sets (or systems) of various kinds, where the method of Buchberger-Gr\"obner bases is
used as a black-box tool to handle some of the involved algebraic computational issues.
The second direction is devoted to the investigation of alternative algorithms from
the constructive theory of partial differential equations developed by
C.\,H.\,Riquier, M.\,Janet, J.\,F.\,Ritt,  J.\,M.\,Thomas, and others for
efficient computation of Gr\"obner bases with variants. The literature is rich for each of these directions and there is a large amount of work which may be considered as relevant to what is presented here. A review of such related work is beyond the intended scope of this paper.

The connection we have established between Ritt's characteristic sets and Buchberger's Gr\"obner bases is expected to stimulate further research and development on some of the outstanding problems in the above-mentioned directions. For example, the pseudo-divisibility relations shown in Theorem~\ref{mainth} (c) and (d) allow us to split the Buchberger-Gr\"obner basis $\pset{G}$ by using the explicit and nontrivial factorizations of $C_{l}$ and $C_{k+1}$ to compute a normal triangular decomposition. Let us explain (part of) the splitting process briefly.

 1. When $I_{k+1}$ is not R-reduced with respect to $C_l$, the first conclusion of Theorem~\ref{mainth} (c) implies that $I_1^{t_1}\cdots I_{l}^{t_l}I_{k+1}\in\langle\pset{C}_l\rangle\subset\langle\pset{G}\rangle$ for some integers $t_i\geq 0$, where $I_i=\ini(C_i)$ for $1\leq i\leq r$. Note that $I_1,\ldots, I_{k+1}$ are all B-reduced with respect to $\pset{G}$. Now $\pset{G}$ can be split into $\pset{G}\cup\{I_1\},\ldots, \pset{G}\cup\{I_l\}, \pset{G}\cup\{I_{k+1}\}$.

2. To deal with the case in which $I_{k+1}$ is R-reduced with respect to $C_l$, we recall
  the pseudo-remainder formula $I^qC_{l}=QI_{k+1}+L$ for $\prem(C_{l}, I_{k+1}, y_{l})=L$, where  $I=\ini(I_{k+1})$ and $q\geq 0$. One can easily see that $\ini(Q)=I^{q-1}I_{l}$. Suppose that $\prem(I^{q-1}I_{l}, \pset{C}_{l-1})=I_{l}\prem(I^{q-1}, \pset{C}_{l-1})\equiv 0$, so $\prem(I^{q-1}, \pset{C}_{l-1})\equiv 0$. Then according to \eqref{rmdr:TS} there exist nonnegative integers $s_1,\ldots,s_{l-1}$ such that $I_1^{s_1}\cdots I_{l-1}^{s_{l-1}}I^{q-1}\in\langle\pset{C}_{l-1}\rangle\subset\langle\pset{G}\rangle$. Of course, $I$ is B-reduced with respect to $\pset{G}$. Now $\pset{G}$ can be split into $\pset{G}\cup\{I_1\},\ldots, \pset{G}\cup\{I_{l-1}\}, \pset{G}\cup\{I\}$.

 3. According to the first conclusion of Theorem~\ref{mainth} (d), $\prem(L, \pset{C}_{l-1})\equiv 0$. This implies that $I^qC_{l}-QI_{k+1}\in\sat(\pset{C}_{l-1})$ and thus $QI_{k+1}\in\sat(\pset{C}_{l})$. Therefore, $I_1^{q_1}\cdots I_{l-1}^{q_{l-1}}QI_{k+1}\in\langle\pset{C}_{l}\rangle\subset\langle\pset{G}\rangle$
 for some integers $q_i\geq 0$. As $\deg(I_{k+1},y_{{l}})>0$, $Q$ is obviously R-reduced with respect to $C_l$.
 Now suppose that $\prem(\ini(Q), \pset{C}_{l-1})=\prem(I^{q-1}I_{l}, \pset{C}_{l-1})\not\equiv 0$.
 It follows from Lemma~\ref{satlem} that $\prem(Q, \pset{C}_{l-1})\not \equiv 0$. Clearly, $\prem(Q, \pset{C}_{l-1})$ is R-reduced with respect to $\pset{C}_{l}$ and
 thus B-reduced with respect to $\pset{G}$. Then $\pset{G}$ can be split into $\pset{G}\cup\{I_1\},\ldots, \pset{G}\cup\{I_l\}, \pset{G}\cup\{\prem(Q, \pset{C}_{l-1})\}, \pset{G}\cup\{I_{k+1}\}$.

In any case of splitting, the split polynomial set $\pset{G}^+$ is obtained from $\pset{G}$ by adjoining a nonzero polynomial $F\in\fk[u_1,\ldots,u_m,y_1,\ldots,y_{l}]$ which is B-reduced with respect to $\pset{G}$ and thus does not belong to $\langle\pset{G}\rangle$. Hence $\langle\pset{G}\rangle\subset \langle\pset{G}^+\rangle$ and $\langle\pset{G}\rangle\neq \langle\pset{G}^+\rangle$.

Consider further the reduced Buchberger-Gr\"obner basis and the W-characteristic set of each $\langle\pset{G}^+\rangle$ and continue the splitting process. In view of the Ascending Chain Condition of polynomial ideals \cite[Ch.\,2, Theorem~7]{clo96}, the splitting process must terminate in a finite number of steps.
Finally, we shall obtain finitely
many reduced plex Buchberger-Gr\"obner bases $\pset{G}_1, \ldots, \pset{G}_e$ such that
$\zero(\pset{P})=\zero(\pset{G}_1)\cup\cdots\cup\zero(\pset{G}_e)$, or
equivalently
$\sqrt{\langle\pset{P}\rangle}
=\sqrt{\langle\pset{G}_1\rangle}\cap\cdots\cap\sqrt{\langle\pset{G}_e\rangle}$,
and
the W-characteristic set $\pset{C}_i$ of each $\langle\pset{G}_i\rangle$ is normal
for $1\leq i\leq e$.
In the case when $\sat(\pset{C}_j)\neq\langle\pset{G}_j\rangle$ which can be determined, for instance,
by computing the reduced plex Buchberger-Gr\"obner basis $\pset{G}_j'$ of $\sat(\pset{C}_j)$,
$\pset{G}_j$ may be further split into $\pset{G}_j'$ (called a \emph{strong} regular Buchberger-Gr\"obner basis)
and $\pset{G}_j\cup\{F_j\}$, where $F_j$ is the
product of the initials of the polynomials in $\pset{C}_j$. Therefore,
we may also ensure that $\sat(\pset{C}_i)=\langle\pset{G}_i\rangle$ for all $i$ ($1\leq i\leq e$).

The process described above can be formulated as an algorithm.
We shall detail the algorithm and discuss computational aspects elsewhere.

We have no idea how to extend the presented connection from the algebraic to the differential
case because a general theory of Gr\"obner bases for differential polynomial ideals is still lacking.
This paper has benefited from the discussions which the author had with Xiaoliang Li, Chenqi Mou, and Jing Yang.

\end{document}